\documentclass[titlepage,amsmath,amssymb,amsfonts]{article}
\usepackage{amsfonts,amssymb,amsmath}
\evensidemargin\oddsidemargin
\newcommand{\PSbox}[3]{\mbox{\rule{0in}{#3}\includegraphics{#1}\hspace{#2}}}

\addtocounter{section}{-1}

\newcommand\Endom{{\mathcal {END}}}
\newcommand\K{{\mathcal K}}
\newcommand\LC{{\mathcal {LC}}}
\newcommand\TB{{\mathcal {TB}}}
\newcommand\Assoc{{\mathcal {ASSOC}}}
\newcommand\Comm{{\mathcal {COMM}}}
\newcommand\BV{{\mathcal {BV}}}
\newcommand\Lie{{\mathcal {LIE}}}
\newcommand\Poiss{{\mathcal {POISS}}}
\newcommand\Gerst{{\mathcal {GERST}}}
\newcommand\g{{\mathfrak g}}
\renewcommand\O{{\mathcal O}}
\newcommand\Rd{{\mathbb R}^d}
\newcommand\R{{\mathbb R}}
\newcommand\Z{{\mathbb Z}}
\newcommand\Q{{\mathbb Q}}
\newcommand\kk{\Bbbk}
\newcommand\BSD{{\mathcal {BSD}}}

\newtheorem{theorem}{Theorem}[section]
\newtheorem{lemma}[theorem]{Lemma}
\newtheorem{conjecture}[theorem]{Conjecture}
\newtheorem{proposition}[theorem]{Proposition}

\newtheorem{definition}[theorem]{Definition}

\newtheorem{corollary}[theorem]{Corollary}

\newcommand\rth{\refstepcounter{theorem}}
\newcommand\numb{\rth\arabic{section}.\arabic{theorem}}

\newtheorem{example}[theorem]{Example}

\newtheorem{remark}[theorem]{Remark}

\begin{document}
\large
\author{Victor Tourtchine\thanks{The autor is partially supported
by the grant RFBR 00-15-96084}}
\title{On the homology of the spaces of long knots}
\date{}
\maketitle

\begin{abstract}


\vspace{1mm}

This paper is a more detailed version of~\cite{T1},
where the first term of the Vassiliev spectral sequence
(computing the homology of the space of long knots in ${\mathbb R}^d$,
$d\ge 3$) was described in terms of  the Hochschild homology of the
Poisson algebras operad for $d$ odd, and of the Gerstenhaber algebras
operad for $d$ even. In particular, the bialgebra of chord diagrams
arises as some subspace of this homology. The homology in question
is the space of characteristic classes for Hochschild cohomology
of Poisson (resp. Gerstenhaber) algebras considered as associative
algebras. The paper begins with necessary preliminaries on operads.

Also we give a simplification of the computations of the first
term of the Vassiliev
spectral sequence.

We do not give proofs of the results.

\end{abstract}

\section{Introduction}\label{s0}

First we recall some known facts on the Vassiliev spectral sequence and
then proceed to explaining of the main idea of the work.

\subsection{}

Let us fix a non-trivial linear map
$l:{\mathbb R}^1\hookrightarrow{\mathbb R}^d$.  We will consider the {\it
space of long knots}, {\it i.~e.}, of injective smooth non-singular
maps ${\mathbb
R}^1\hookrightarrow{\mathbb R}^d$,
that coincide with the map $l$ outside
some compact set (this set is not fixed). The long knots form an open
everywhere dense subset in the affine space ${\mathcal K}$ of all smooth maps
${\mathbb R}^1\to {\mathbb R}^d$ with the same behavior at infinity.
The complement $\Sigma\subset{\mathcal K}$ of this dense subset is called
the {\it discriminant space}. It
consists of the maps having self-intersections or singularities.
Any cohomology class $\gamma\in H^i({\mathcal K}\backslash\Sigma)$
of the knot space can be realized as the linking coefficient with an
appropriate chain in $\Sigma$ of codimension $i+1$ in ${\mathcal K}$.

Following~\cite{V5} we will assume that the space ${\mathcal K}$ has a
very large but finite dimension $\omega$. A partial justification of this
assumption uses finite dimensional approximations of $\mathcal K$,
see~\cite{V1}. 
Below we indicate by quotes non-rigorous assertions using this
assumption and needing a reference to \cite{V2} for such a justification.

The main tool of Vassiliev's approach to computation of the (co)homology
of the
knot space is the {\it simplicial resolution} $\sigma$ (constructed in
\cite{V1})
of the discriminant $\Sigma$. This resolution is also called the {\it
resolved discriminant}.
The natural projection $\Pi:\bar\sigma\to\bar\Sigma$ is a ``homotopy
equivalence'' between the ``one-point compactifications'' of the spaces
$\sigma$ and
$\Sigma$. By the ``Alexander duality'', the reduced homology groups
$\tilde H_*(\bar\sigma,\kk)\equiv\tilde H_*(\bar\Sigma,\kk)$ of
these compactifications ``coincide'' (up to a change of dimension) with the
cohomology groups of the space of knots:

$$
\tilde H^i({\mathcal K}\backslash\Sigma,\kk)\simeq
\tilde H_{\omega -i-1}(\bar\Sigma,\kk) 
\equiv\tilde H_{\omega-i-1}(\bar\sigma,\kk),
\eqno(\numb)\label{eq01}
$$
where $\kk$ is a commutative ring of coefficients.

In the space $\sigma$ there is a natural filtration
$$
\emptyset=\sigma_0\subset\sigma_1\subset\sigma_2\subset\dots.
\eqno(\numb)\label{eq02}
$$

\medskip

\begin{conjecture}\label{con03}
The spectral sequence (called {\sl Vassiliev's main
spectral sequence}) associated with the filtration~(\ref{eq02}) and computing
the
``Borel-Moore homology groups of the resolution $\sigma$''
stabilizes over $\mathbb Q$ in
the first
term.~$\square$
\end{conjecture}

\begin{conjecture}\label{con04}
{\rm (due to Vassiliev)}
Filtration~(\ref{eq02}) ``homotopically
splits'',
in other words, $\bar\sigma$ is ``homotopy equivalent'' to the wedge
$\bigvee _{i=1}^{+\infty}(\bar\sigma_i/\bar\sigma_{i-1})$.~$\square$
\end{conjecture}

This conjecture would imply the stabilization 
of our {\it main spectral
sequence} in the first term over any commutative ring $\kk$ of coefficients.

Due to the Alexander duality, the filtration~(\ref{eq02})
induces the filtrations

$$
H_{(0)}^*({\mathcal K}\backslash\Sigma)\subset H_{(1)}^*({\mathcal
K}\backslash\Sigma)\subset
H_{(2)}^*({\mathcal K}\backslash\Sigma)\subset\dots,
\eqno(\numb)\label{eq05}
$$

$$
H_*^{(0)}({\mathcal K}\backslash\Sigma)\supset H_*^{(1)}({\mathcal
K}\backslash\Sigma)\supset
H_*^{(2)}({\mathcal K}\backslash\Sigma)\supset\dots
\eqno(\numb)\label{eq06}
$$
in respectively the cohomology and homology  groups of the space of knots. For
$d\ge 4$
the
filtrations~\eqref{eq05},~\eqref{eq06} are finite for any dimension $*$.
The Vassiliev spectral sequence
in this case
computes the graded quotient associated with the above filtrations.

In the most intriguing case $d=3$ almost nothing is clear. For the dimension
$*=0$ the
filtration~(\ref{eq05}) does not exhaust the whole cohomology of degree zero.
The knot invariants obtained by this method are called the
Vassiliev invariants, or invariants
of finite type. One can define them in
a more simple and geometrical way, see~\cite{V1}.
The dual space to the graded quotient of the space of finite type knot
invariants
is
the {\it bialgebra of chord diagrams}. The invariants and the bialgebra in
question
were intensively studied in the last decade, see~\cite{AF, BN1, CCL, ChD,
GMM, HV, K, L,
P, S, Vai, Z}.
The completeness conjecture for the Vassiliev knot invariants is the question
about the convergence of the filtration~(\ref{eq06}) to zero  
for $d=3$, $*=0$. The
realization theorem
of M.~Kontsevich \cite{K} proves that the Vassiliev spectral
sequence for $d=3$,
$*=0$
also computes the corresponding associated quotient (for positive dimensions $*$
in the case
$d=3$ even this is not for sure) and does stabilize in the first term. 
The groups
of the
associated graded quotient to filtrations~(\ref{eq05}),~(\ref{eq06}) in the case $d=3$, $*>0$
are
some quotient groups
of the groups calculated by  Vassiliev's main spectral sequence.

Let us mention that there is a natural way to construct real
cohomology
classes of the knot spaces by means of
configuration space integrals, that generalizes
in a non-trivial way the Vassiliev knot invariants obtained in three
dimensions from the Chern-Simons perturbation theory, see~\cite{AF, BN2,
CCL, GMM, P}.
May be this approach leads to a proof of Conjecture~\ref{con03}.

To compute the first term of the main spectral sequence, V.~A.~Vassiliev
introduced an {\it auxiliary filtration} in the spaces $\sigma_i\backslash
\sigma_{i-1}$, see~\cite{V1, V4}. The {\it auxiliary spectral
sequence} associated to this filtration degenerates in the first term, because
its first term (for any $i$) is concentrated at only one line.
Therefore {\it the second term of the auxiliary
spectral sequence is isomorphic
to the first term of the main spectral sequence}. The term $E_0^{*,*}$ of the
auxiliary spectral sequence together with its differential (of degree zero)
is a direct sum of tensor products of complexes of connected graphs.
The homology groups of the complex of connected graphs with $m$ labelled
vertices are  concentrated in the dimension $(m-2)$ only, and 
the only non-trivial group is isomorphic
to ${\mathbb Z}^{(m-1)!}$, see~\cite{V3, V4}.  This homology group has a
nice description as the quotient by the 3-term relations of the space spanned
by trees (with $m$ labelled vertices), see~\cite{T2}.

\subsection{}\label{s02}

In fact V.~A.~Vassiliev considered only \underline{cohomological} case
of the main and auxiliary spectral sequences, {\it i.~e.}, the case
corresponding
to the homology of the discriminant and (by the Alexander duality) to the
\underline{cohomology} of the knot space ${\mathcal K}\backslash\Sigma$.
This was because a convenient description only of the homology of complexes
of connected graphs was known.
It was noticed in~\cite{T3} that
the cohomology of complexes
of connected graphs (with $m$ vertices) admits also a very nice description as
the $m$-th component of the Lie algebras operad. This isomorphism comes from
the following observation.

Let us consider the space $Inj(M,{\mathbb R}^d)$ of injective maps from a finite
set $M$ of cardinality $m$ into ${\mathbb R}^d$, $d\ge 2$. This space can be
viewed
as a fi\-ni\-te-\-di\-men\-sio\-nal analogue of the knot spaces.
The corresponding discriminant (consisting of non-injective maps $M\to\Rd$)
has also a simplicial resolution, whose filtration (analogous to~(\ref{eq02}))
does split homotopically, see~\cite{V2, V4}. The superior non-trivial
term $\sigma_{m-1}\backslash\sigma_{m-2}$ of the filtration provides exactly
the complex of connected graphs with $m$ vertices labelled by the elements
of $M$.  Its only non-trivial homology group is isomorphic by 
the Alexander duality to the cohomology group
in the maximal degree of the space $Inj(M,\Rd)$. On the other hand the above
space is homotopy equivalent to the $m$-th space of the little cubes operads,
see Section~\ref{s4}. The homology operad of this topological
operad is well known, see~\cite{Co}
(and Section~\ref{s4}). For different $d$ of the same parity this homology 
operad is
the same up to a change of grading and for odd (resp. even) $d$ it 
is the
Poisson (resp. Gerstenhaber) algebras operad containing (in the maximal
degree)
the operad of Lie algebras.

An analogous periodicity takes place for the spaces of knots. The degree zero
term of the main spectral sequence together with its differential
(of degree zero) depends up to a change of grading on the parity
of $d$ only. Obviously the same is true for the whole auxiliary spectral
sequence.

The above description of the cohomology of  complexes of connected graphs
allows one to describe easily the Vassiliev spectral sequence in the
\underline{homological} case.
The main results of these computations are explained
in Section~\ref{s5}
(see Theorems~\ref{t54},~\ref{t58},~\ref{t59},~\ref{t510}).
The proofs in full detail will be
given in~\cite{T2}, (see also the
Russian  version~\cite{T3}).

\subsection{}

The paper is organized as follows.

In Sections~\ref{s1},~\ref{s2},~\ref{s3}
we give some preliminaries
on  linear graded operads. We give a short definition and examples
(that will be  useful  for us) of linear graded operads
(Section~\ref{s1}).  We construct a graded Lie algebra structure on the space
of any linear graded operad (Section~\ref{s2}). We define
a Hochschild complex
for any graded linear operad supplied with a morphism from the associative
algebras operad to this operad (Section~\ref{s3}). The content of
Sections~\ref{s2} and~\ref{s3} was borrowed (up to a slightly different
definition of signs) from~\cite{GV}.  

Section~\ref{s4} is devoted to the May
operad of little cubes and to its homology. As it was already mentioned the
m-th component of this operad is a space homotopy equvalent to the space
$Inj(\{1,2,\dots,m\},\Rd)$ of injective maps $\{1,2,\dots,m\}
\hookrightarrow\Rd$.
In Section~\ref{s4}  we explain how the stratification in the
discriminant set of non-injective maps $\{1,2,\dots,m\}\to\Rd$
corresponds to a direct sum decomposition of the homology
$H_*\bigl(Inj(\{1,2,\dots,m\},\Rd)\bigr)$.

In Section~\ref{s5} we describe a natural
stratification in the discriminant set
$\Sigma$ of long singular ``knots''. This stratification provides a direct
sum decomposition of the first term of the
auxiliary spectral sequence. In the case
of even $d$ the first term of Vassiliev's auxiliary spectral sequence
is completely described by the following theorem:

\medskip

{\bf Theorem~\ref{t54}.} {\it The first term of Vassiliev's
homological auxiliary
spectral sequence together with its first differential is isomorphic
to the normalized Hochschild complex of the Batalin-Vilkovisky algebras
operad.}~$\square$
\medskip

Unfortunately in the case of odd $d$ it is not possible to describe the
corresponding complex in terms of the  Hochschild complex for some graded
linear operad. A description of this complex is given 
in Section~\ref{s8}.
Nevertheless the
homology of the complex in question ({\it i.~e.}, the first
term of the main spectral
sequence) over $\Q$ can be defined in terms of the 
Hochschild homology of the
Poisson algebras operad in the case of odd $d$ (and of the Gerstenhaber
algebras operad in the case of even $d$). A precise statement is given by
Theorem~\ref{t510}.

In Section~\ref{s5} we also introduce complexes homologically equivalent
(for any
commutative ring $\kk$ of coefficients) to the first term of the auxiliary
spectral sequence. These complexes simplify a lot the computations of the
second term.

In Section~\ref{s6} we explain how the bialgebra of chord diagrams
arises in our construction. We formulate some problems concerning it.

Section~\ref{s7} does not contain any new results and serves rather to explain
one remark of M.~Kontsevich. We study there the homology operads of some
topological operads, that we call {\it operads of turning balls}. 
These homology operads
make more clear the difference that we have in the cases of odd and even $d$.

In Section~\ref{s8} we describe the first term of the auxiliary spectral
sequence together with the degree $1$ differential both for $d$ even and odd.
The corresponding complex is called {\it Complex of bracket star-diagrams}.

In Section~\ref{s9} we construct a differential bialgebra structure on
this complex. We conjecture that this differential bialgebra structure is 
compatible with the homology bialgebra structure of the space of long knots.

\bigskip

\underline{\bf Acknowledgements.} I would like to thank my scientific advisor
V.~Vassiliev for his consultations and for his support during the work.
Also I would like to express my gratitude to Ecole Normale Sup\'erieur
and  Institut des Hautes Etudes Scientifiques for hospitality.
I would like to thank M.~Kontsevich for his attention to this work.

Also 
I am  grateful to P.~Cartier, A.~V.~Chernavsky, M.~Deza,
D.~Panov, M.~Finkelberg, S.~Loktev, I.~Marin, G.~Racinet, A.~Stoyanovsky.

\section{Linear operads}\label{s1}
The definition of many algebraic structures  on a linear space (such as the 
commutative, associative, Lie algebra structures) consists of setting
several polylinear operations (in these three cases, only one binary 
operation), 
that should satisfy some composition identities 
(in our example, associativity or Jacoby identity). 
Instead of doing this one can
consider the spaces of all polylinear $n$-ary operations, for all $n\ge 0$,
and the composition rules, that arize from the corresponding algebraic
structure. The natural formalization of this object is given by the notion 
of operad.

\medskip
\centerline{\underline{\bf Definition}}
\medskip

Let $\kk$ be a commutative ring of coefficients. A {\it graded $\kk$-linear
operad}
$\O$ is a collection $\{\O(n),\, n\ge 0\}$ of graded
$\kk$-vector spaces
equipped with the following set of data:

\medskip
(i) An action of the symmtric group $S_n$ on $\O(n)$ for each $n\ge 2$.

(ii) Linear maps (called {\it compositions}), preserving the grading,
$$
\gamma_{m_1,\dots,m_l}:\O(l)\otimes\bigl(\O(m_1)\otimes\dots\otimes\O(m_l)
\bigr)
\to \O(m_1+\dots +m_l)
\eqno(\numb)\label{eq11}
$$
for all $m_1,\dots,m_l\ge 0$. We write $\mu(\nu_1,\dots,\nu_l)$ instead of
$\gamma_{m_1,\dots,m_l}(\mu\otimes\nu_1\dots\otimes\nu_l)$.

(iii) An element $id\in\O(1)$, called the {\it unit}, such that $id(\mu)=
\mu(id,\dots,id)=\mu$ for any non-negative $l$ and any $\mu\in\O(l)$.

\medskip
It is required that these data satisfy some conditions of associativity and
equivariance with respect to the symmetric group actions,
see~\cite[Chapter~1]{May},~\cite{HS},~\cite{Lo}.

One can consider any element $x\in\O(l)$ as something that has $l$ inputs
and 1 output:

\hspace*{67mm}
\PSbox{tree.pstex}{15mm}{48mm}
\vspace{3mm}\vspace{1mm}

\nopagebreak

\parbox[b]{163mm}{
\centerline{(Figure \numb)}\label{f12}
}
\vspace{3mm}

The composition operation~(\ref{eq11}) of $x\in\O(l)$ with $y_1\in\O(m_1)$,
$\dots$, $y_l\in\O(m_l)$ is the substitution of $y_1$, $\dots$,
$y_l$ into $l$ inputs of $x$:

\hspace*{60mm}
\PSbox{compos.pstex}{15mm}{70mm}
\LARGE
\begin{picture}(50,10)
\put(-130,85){$x(y_1,y_2,\dots,y_l)$}
\end{picture}
\large
\vspace{3mm}\vspace{1mm}

\nopagebreak

\parbox[b]{163mm}{
\centerline{(Figure \numb)}\label{f13}
}
\vspace{3mm}

The resulting element has $m_1+\dots +m_l$ inputs and 1 output, {\it i.~e.},
it belongs to $\O(m_1+\dots +m_l)$.

\medskip

\centerline{\underline{\bf Examples}}

\medskip

We will give several well known examples of graded linear operads, 
see~\cite{GK,Co,G}.

\vspace{2mm}

{\Large\bf a)} Let $V$ be a graded $\kk$-vector space.
We define the {\it endomorphism
operad} $\Endom(V):=\{Hom(V^{\otimes n},V),n\ge 0\}$. The unit element
$id\in Hom(V,V)$ is put to be the identical map $V\to V$.
The composition operations and the symmetric group actions are defined in
the canonical way.

This operad is very important due to the following definition.

\begin{definition}
\label{d14}
{\rm
Let $\O$ be a graded linear operad. By an
{\it $\O$-algebra} (or {\it algebra over $\O$}) 
we call any couple $(V,\rho)$, where $V$ is a graded vector space
and $\rho$ is  a morphism $\rho:\O\to \Endom(V)$ of operads.}~$\square$
\end{definition}

In other words, the theory of algebras over an operad $\O$
is the representation theory of $\O$.

\vspace{3mm}

{\Large\bf b)} The {\it operad $\Lie$ of Lie algebras}.
The component $\Lie(0)$ of this
operad is trivial.
The $n$-th component $\Lie(n)$
is defined as the subspace of the free Lie algebra
with  generators $x_1,x_2,\dots,x_n$,  that is spanned
 by the brackets containing each generator exactly once.

\begin{example}\label{ex15}
{\rm
For $n=5$ one can take the bracket
$[[[x_5,x_3],[x_1,x_2]],x_4]$ as an element of $\Lie(5)$.~$\square$
}
\end{example}

Since we work in the category of graded vector spaces, we need
to define a grading on each considered space. 
The grading of the spaces $\Lie(n)$, $n\ge 1$, is put to be zero.
It is well known that these spaces are free $\kk$-modules,
$\Lie(n)\simeq\kk^{(n-1)!}$. The $S_n$-action is defined by
permutations of $x_1,\dots,x_n$.

Let $A(x_1,\dots,x_l)$, $B_1(x_1,\dots,x_{m_1})$, $\dots$,
$B_l(x_1,\dots,x_{m_l})$ be brackets respectively from
$\Lie(l)$, $\Lie(m_1)$, $\dots$, $\Lie(m_l)$. We define the composition
operations~(\ref{eq11}) as follows.
$$
A(B_1,\dots,B_l)(x_1,\dots,x_{m_1+\dots +m_l}):=A\bigl(
B_1(x_1,\dots,x_{m_1}),B_2(x_{m_1+1},\dots,x_{m_1+m_2}),\dots,
$$
$$
B_l(x_{m_1+\dots+m_{l-1}+1},\dots,x_{m_1+\dots +m_l})\bigr).
\eqno(\numb)\label{eq16}
$$

The element $x_1\in\Lie(1)$ is the unit element for
this operad.

Note that a $\Lie$-algebra structure in the sense of Definition~\ref{d14}
is exactly the same as a (graded) Lie algebra structure in the usual sense.
Indeed, the element $\rho([x_1,x_2])\in Hom(V^{\otimes 2},V)$ always
defines a Lie bracket. The converse is also true. It is easy to see
that the element $[x_1,x_2]\in\Lie(2)$ generates
the operad $\Lie$, so if we put $\rho([x_1,x_2])$ equal to our
Lie bracket, then we immediately obtain a map $\rho:\Lie\to \Endom(V)$ of the
whole operad $\Lie$.

\vspace{3mm}

{\Large\bf c)} By analogy with the operad $\Lie$
one defines the {\it operad} $\Comm$ (resp. $\Assoc$) {\it of commutative} 
(resp.
{\it associative}) {\it algebras}. The space
$\Comm(n)$ (resp. $\Assoc(n)$), for $n\ge 1$, is put to be
one-dimensionl (resp. $n!$-dimensional) free $\kk$-module
defined as the subspace of the free commutative (resp. associative)
algebra with ge\-ne\-ra\-tors $x_1,\dots,x_n$, that is spanned
by the monomials containing each generator exactly once.
The $S_n$-actions and the composition operations are defined in the
same way as for the operad $\Lie$. The element $x_1\in\Comm(1)$
(resp. $\Assoc(1)$) is the unit element $id$.
There are two different ways to
define the space $\Comm(0)$ (resp. $\Assoc(0)$). We can put this space
to be trivial or one-dimensional. In the first case we get the
operad of commutative (resp. associative) algebras without unit,
in the second case -- the operad of commutative (resp. associative)
algebras with unit. Below we will consider the second situation.

\begin{remark}\label{r17}
{\rm
A commutative (resp. associative) algebra structure
is exactly the same as a structure of $\Comm$-algebra (resp. $\Assoc$-algebra)
in the sense of Definition~\ref{d14}.~$\square$
}
\end{remark}

\vspace{3mm}

{\Large\bf d)} The  {\it operads} $\Poiss$, $\Gerst$, $\Poiss_d$
{\it of Poisson, Gerstenhaber} and {\it $d$-Poisson algebras}.
First of all, let us recall the definition of
Poisson, Gerstenhaber and
$d$-Poisson algebras.

\begin{definition}\label{d18}
{\rm
A graded commutative algebra $A$ is called a
{\it $d$-Poisson algebra}, if it has a Lie bracket
$$
[.,.]:A\otimes A\to A
$$
of degree $-d$. The bracket is supposed to be compatible
with the multiplication. This means that for any elements
$x,y,z\in A$

$$
[x,yz]=[x,y]z+(-1)^{\tilde{y}(\tilde{x}-d)}y[x,z].~\square
\eqno(\numb)\label{eq19}
$$
}
\end{definition}

0-Poisson (resp. 1-Poisson) algebras are called
simply {\it Poisson} (resp. {\it Gerstenhaber}) algebras.
1-Poisson algebras are called Gerstenhaber algebras in honor of
Murray Gers\-ten\-haber, who discovered this structure on the
Hochschild cohomology of associative algebras, see~\cite{Ge} and also
Sections~\ref{s2} and \ref{s3}.

\begin{example}\label{ex110}
{\rm
Let $\g$ be a graded Lie algebra with the bracket
of degree $-d$. Then the symmetric (in the graded sense) algebra
$S^*\g$ has a natural structure of a $d$-Poisson algebra with
the usual multiplication of a symmetric algebra and with the bracket
defined by the following formula:
$$
[A_1\cdot A_2\dots A_k,B_1\cdot B_2\dots B_l]=\sum_{i,j}(-1)^
\epsilon A_1\dots\widehat{A_i}\dots A_k\cdot [A_i,B_j]\cdot B_1\dots
\widehat{B_j}\dots B_l,
\eqno(\numb)\label{eq111}
$$
where $\epsilon=\tilde A_i(\sum_{i'=i+1}^k \tilde A_{i'}) +
\tilde B_j(\sum_{j'=1}^{j-1}\tilde B_{j'}).$~$\square$
}
\end{example}

\begin{remark}\label{r112}
{\rm
For a graded Lie algebra $\g$ with the
bracket of degree 0 the $d$-tuple suspension $\g[d]$ is also a
graded Lie algebra with the bracket of degree $-d$. Thus, the space
$S^*(\g[d])$ is a $d$-Poisson algebra.~$\square$
}
\end{remark}

Note, that for commutative, associative, or Lie algebras the $n$-th
component of the ope\-rad is defined as the space of all natural
polylinear $n$-ary operations, that come from the corresponding
algebra structure. Now, let us describe the spaces $\Poiss_d(n)$
of all natural polylinear $n$-ary operations of $d$-Poisson algebras.
Consider a free graded Lie algebra $Lie_d(x_1,\dots,x_n)$
with the bracket of degree $-d$ and with the generators $x_1,\dots,
x_n$ of degree zero, and consider
(following Example~\ref{ex110}) the $d$-Poisson algebra
$Poiss_d(x_1,\dots,
x_n):=
S^*\bigl( Lie_d(x_1,
\dots,x_n)\bigr)$. This is a free $d$-Poisson algebra. We will
define the space $\Poiss_d(n)$ as the subspace of $Poiss_d(x_1,\dots,
x_n)$ spanned by the products (of brackets) containing each generator
$x_i$ exactly once. For instance for $n=5$ we will take the pro\-duct
$[x_1,x_3]\cdot[[x_2,x_5]x_4]$ as an element of $\Poiss_d(5)$.
The unit element $id$ is $x_1\in\Poiss_{d}(1)$. The symmetric group actions and
the composition operations are defined analogously to the case of the
operad $\Lie$.

The space $\Poiss_d(n)$ can be decomposed into a direct sum with the summands
numbered by
partitions
of the set $\{1,\dots,n\}$:
$$
\Poiss_d(n)=\bigoplus_A\Poiss_d(A,n).
\eqno(\numb)\label{eq113}
$$
For a partition $A=\{\bar A_1,\dots,\bar A_{\#A}\}$
of the set $\{1,\dots,n\}=\coprod_{i=1}^{\# A}\bar A_i$ we define the
space $\Poiss_d(A,n)\subset\Poiss_d(n)$ to be linearly
spanned by products of
$\#A$ brackets, such that the $i$-th bracket contains generators only
from the set $\bar A_i$ (thus, each generator from $\bar A_i$
is presented exactly once in the $i$-th bracket).

Let $\bar A_1,\dots\bar A_{\#A}$ be of cardinalities $a_1,\dots,
a_{\#A}$ respectively, then
$$
\Poiss_d(A,n)\simeq\otimes_{i=1}^{\#A}\kk^{(a_i-1)!}.
\eqno(\numb)\label{eq114}
$$
This implies  that the space $\Poiss_d(n)$ is isomorphic to
$\kk^{n!}$, and its Poincar\'e polynomial is $(1+t^{-d})(1+2t^{-d})\dots
(1+(n-1)t^{-d}).$

\vspace{3mm}

{\Large\bf e)} The  operad $\BV$ (resp. $\BV_d$, $d$ being odd) of
Batalin-Vilkovisky (resp. $d$-Batalin-Vilkovisky) algebras.
\begin{definition}\label{d115}
{\rm
A Gerstenhaber algebra (resp. $d$-Poisson
algebra, for odd $d$) $A$ is called a {\it Batalin-Vilkovisky algebra}
(resp. {\it $d$-Batalin-Vilkovisky algebra}), if $A$ is supplied with
a linear map $\delta$ of degree -1 (resp. $-d$)
$$
\delta:A\to A,
$$
such that

(i) $\delta^2=0$,

(ii) $\delta(ab)=\delta(a)b+(-1)^{\tilde a}a\delta(b)+(-1)^{\tilde a}
[a,b].$~$\square$
}
\end{definition}

Note that (i) and (ii) imply

(iii) $\delta([a,b])=[\delta(a),b]+(-1)^{\tilde a +1}[a,\delta(b)].$

\begin{example}\label{ex116}
{\rm
Let $\g$ be a graded Lie algebra with the bracket of
degree zero. Then the exterior algebra $\Lambda^*\g:=S^*(\g[1])$ is
a Batalin-Vilkovisky algebra, where the structure of a Gerstenhaber
algebra is from Remark~\ref{r112}; the operator $\delta$ is the
standard differential on the chain-complex $\Lambda^*\g$:
$$
\delta(A_1\wedge\dots\wedge A_k)=\sum_{i<j}(-1)^\epsilon[A_i,A_j]\wedge
A_1\dots \widehat{A_i}\dots\widehat{A_j}\dots\wedge A_k,
\eqno(\numb)\label{eq117}
$$
where $A_1,\dots,A_k\in\g$,
$\epsilon=\tilde A_i+(\tilde A_i+1)(\tilde A_1+\dots
+\tilde A_{i-1}+i-1)+(\tilde A_j+1)(\tilde A_1+\dots+\widehat{\tilde A_i}+\dots
+\tilde A_{j-1} +j-2)$.
In the same way $\g$ defines  the $d$-Batalin-Vilkovisky algebra
$S^*(\g[d])$ for any
odd  $d$.~$\square$
}
\end{example}

Let us describe the $n$-th component of the corresponding operad,
denoted by $\BV$ and $\BV_d$ ($d$ being always odd). Since
the spaces $\BV(n)$ and $\BV_d(n)$ are isomorphic (in the super-sense)
up to a change
of grading, we will consider now only the case $d=1$.
Obviously, the space $\BV(n)$ of all natural polylinear
$n$-ary operations for such algebras contains $\Gerst(n)$.
Consider the symmetric algebra of the free graded Lie algebra
$Lie_1(x_1,\dots,x_n,\delta(x_1),\dots,\delta(x_n))$ with the bracket
of degree -1 and with the generators $x_1,\dots,x_n$ of degree zero
and the generators $\delta(x_1),\dots,\delta(x_n)$ of degree -1.
This space has a structure of a Batalin-Vilkovisky algebra. (In fact
it is a free $\BV$-algebra with generators $x_1, \dots,x_n$). In
$S^*\bigl( Lie_1(x_1,\dots,x_n,\delta(x_1),\dots,\delta(x_n))\bigr)$
we will take the subspace $\BV(n)$ linearly spanned by all the
products (of brackets), containing each index $i\in\{1,\dots,n\}$
exactly once. For instance $[\delta(x_1),x_3]\cdot[x_2,[\delta(x_4),
\delta(x_5)]]$ belongs to $\BV(5)$. Due to relations (i), (ii), (iii)
for $\delta$, this subspace is exactly the space of all natural polylinear
$n$-ary operations on
Batalin-Vilkovisly algebras.

Analogously to~(\ref{eq113}) the space $\BV(n)$ can be
decomposed into  the direct sum:
$$
\BV(n)=\bigoplus_{A,S}\BV(A,S,n),
\eqno(\numb)\label{eq118}
$$
where
$A$ is a partition of the set $\{1,\dots,n\}$, $S$ is a subset of
$\{1,\dots,n\}$ corresponding to indices $i$ presented by
$\delta(x_i)$ in products of brackets of $\BV(n)$.

Note that the space $\BV(n)$ is isomorphic to $\kk^{2^nn!}$.
The Poincar\'e polynomial of this graded space is
$(1+t^{-d})^{n+1}(1+2t^{-d})(1+3t^{-d})\dots (1+(n-1)t^{-d})$.

\vspace{1.5mm}

In the sequel we will use the following definition.

\begin{definition}\label{d119}
{\rm
For a finite set $M$
any pair $(A,S)$ consisting of a partition $A$ of $M$
and of a subset $S\subset
M$ will be called a {\it star-partition} of the set $M$.
(Each point $i\in S$ will be called  a {\it star}.)~$\square$
}
\end{definition}

Note that we have the following natural morphisms of operads:
$$
\Assoc\to\Comm\to\Poiss_d\dashrightarrow\BV_d
\eqno(\numb\label{eq120})
$$
The last arrow of~\eqref{eq120} is defined only if $d$ is odd.

\section{Graded Lie algebra structure on graded linear
operads}\label{s2}

In this section we define  a graded Lie algebra structure on an arbitrary
graded
li\-near operad. In the next section we use this structure and define
the Hochschild comp\-lex
 for a graded linear operad equipped with a morphism from the
operad $\Assoc$ to our operad. Both these constructions,
which generalize the Hochschild cochain
complex for associative algebras, were introduced in~\cite{GV}.
The only difference between the
operations~\eqref{eq21},~\eqref{bracket},~\eqref{eq32}  given below and those
of~\cite{GV} is in signs. First of all, in the paper~\cite{GV} M.~Gerstenhaber
and A.~Voronov considered only linear (non-graded) ope\-rads,
hence our case
is  more general. But even for the case of purely even gradings
the signs are slightly different. The difference can be easily
obtained by conjugation of the
operations~\eqref{eq21},~\eqref{bracket},~\eqref{eq32} by means of the linear
operator that maps any element $x\in\O(n)$,
$n\ge 0$, to $(-1)^{\frac{n(n-1)}{2}}x$.
\medskip

Let $\O=\{\O(n),n\ge 0\}$ be a graded $\kk$-linear operad. By abuse
of the language the space $\bigoplus_{n\ge 0}\O(n)$ will be also
denoted by $\O$. A tilde over an element will always
designate its grading. For any element $x\in\O(n)$ we put $n_x:=n-1$.
The numbers $n$ and
1 here correspond to $n$ inputs and to 1 output respectively.

Define a new grading $|\,.\,|$ on the space $\O$. For an element
$x\in\O(n)$ we put $|x|:=\tilde x+n_x=\tilde x +n-1$.
It turns out that
$\O$ is a graded Lie algebra with respect to the grading $|\,.\,|$.
Note that the composition operations~(\ref{eq11}) respect this grading.
Define the following collection of multilinear operations on
the space $\O$.
$$
x\{x_1,\dots,x_n\}:=\sum(-1)^\epsilon x(id,\dots,id,x_1,id,\dots,id,
x_n,id,\dots,id)
\eqno(\numb)\label{eq21}
$$
for $x,x_1,\dots,x_n\in\O$, where the summation runs over all possible
substitutions of $x_1,\dots,x_n$ into $x$ in the prescribed order,
$\epsilon:={\sum_{p=1}^nn_{x_p}r_p}+{n_x\sum_{p=1}^n\tilde x_p}+
{\sum_{p<q}n_{x_p}\tilde x_q}$, $r_p$ being the total
number of inputs in $x$ going after $x_p$. For instance, for
$x\in \O(2)$ and arbitrary $x_1,x_2\in\O$
$$
x\{x_1,x_2\}=(-1)^{n_{x_1}+(\tilde x_1+\tilde x_2)+n_{x_1}\tilde x_2}
x(x_1,x_2).
$$

We will also adopt 
 the following convention:
$$
x\{\}:=x.
$$

One can check immediately the following identities:
$$
x\{x_1,\dots,x_m\}\{y_1,\dots,y_n\}=
$$
$$
\sum_{0\le i_1\le j_1\le \dots
\le i_m\le j_m\le n}(-1)^{\epsilon}
x\{y_1,\dots,y_{i_1},x_1\{y_{i_1+1},\dots,y_{j_1}
\},y_{j_1+1},\dots,y_{i_m},
$$
\nopagebreak
$$
x_m\{y_{i_m+1},\dots,y_{j_m}\},y_{j_m+1},\dots,y_n\},
\eqno(\numb)\label{eq22}
$$
where $\epsilon=\sum_{p=1}^m\bigl( |x_p|\sum_{q=1}^{i_p}|y_q|\bigr)$.
(These signs are the same as in~\cite{GV}).

Define a bilinear operation  (respecting the grading 
$|\, .\,|$) $\circ$ on the space $\O$:
$$
x\circ y:=x\{y\},
\eqno(\numb\label{eq22'})
$$
for $x,y\in\O$.

\begin{definition}\label{d22''}
{\rm A graded vector space $A$ with a bilinear operation
$$
\circ:A\otimes A\to A
$$
is called a {\it Pre-Lie algebra}, if for any $x,y,z\in A$ the following holds:
$$
(x\circ y)\circ z -x\circ (y\circ z)=(-1)^{|y||z|}((x\circ z)\circ y-
x\circ (z\circ y))\,\,\square
$$
}
\end{definition}

Any graded Pre-Lie algebra $A$ can be considered as a graded
Lie algebra with the bracket
$$
[x,y]:=x\circ y-(-1)^{|x||y|}y\circ x.
\eqno(\numb\label{bracket})
$$
The description of the operad of Pre-Lie algebras is given in~\cite{Cha}.

The following lemma is a corollary of the identity~(\ref{eq22}) applied
to the case $m=n=1$.

\begin{lemma}\label{l23}
The operation~\eqref{eq22'} defines
a graded Pre-Lie algebra structure on the space
$\O$.~$\square$
\end{lemma}

In particular this lemma means that any graded linear operad $\O$ 
can be considered as a graded Lie algebra with the bracket~\eqref{bracket}.

\section{Hochschild complexes}\label{s3}

Let $\O =\bigoplus_{n\ge0}\O(n)$ be a graded linear operad equipped
with a morphism
$$
\Pi:\Assoc\to\O
$$
from the operad $\Assoc$. This morphism defines the element $m=\Pi
(m_2)\in\O(2)$, where the element $m_2=x_1x_2\in\Assoc(2)$ is
the operation of multiplication. Note that the elements $m_2$, $m$
are odd with respect to the new grading $|\,.\,|$ ($|m|=|m_2|=1$)
and $[m,m]=[\Pi(m_2),\Pi(m_2)]=2\Pi(m_2\circ m_2)=0$. Thus
$\O$ becomes a differential graded Lie algebra with the differential
$\partial$:
$$
\partial x:=[m,x]=m\circ x-(-1)^{|x|}x\circ m,
\eqno(\numb)\label{eq31}
$$
for $x\in \O$.

The complex $(\O,\partial)$ is called the {\it Hochschild complex}
for the operad $\O$. 

\begin{example}\label{ex31'}
{\rm
If $\O$ is the endomorphism operad $\Endom(A)$
of a vector space $A$, and we have a morphism
$$\Pi:\Assoc\to \Endom(A),
$$
that defines an associative algebra structure on $A$, then the
corresponding complex $\bigl(\bigoplus_{n=0}^{+\infty} Hom(A^{\otimes n},A),
\partial\bigr)$ is the standard Hochschild cochain complex of the associative
algebra $A$.~$\square$
}
\end{example}

\begin{example}\label{ex31''}
{\rm
Due to the morphisms~\eqref{eq120} we have the Hochschild complexes
$(\Assoc,\partial)$, $(\Comm,\partial)$, $(\Poiss_d,\partial)$, $(\BV_d,
\partial)$. It can be shown that the complexes $(\Assoc,\partial)$
and $(\Comm,\partial)$ are acyclic.~$\square$
}
\end{example}

Define another grading
$$
deg:=|\,.\,|+1
\eqno(\numb)\label{eq33}
$$
on the space $\O$. With respect
to this grading the bracket $[.,.]$ is homogeneous of degree
$-1$.

It is easy to see that the product $*$, defined as follows
$$
x*y:=(-1)^{|x|}m\{x,y\}=m(x,y),
\eqno(\numb)\label{eq32}
$$
for $x,y\in\O$, together with the differential $\partial$ defines
a differential graded associative algebra structure on $A$ with
respect to the grading $deg=|\,.\,|+1.$

\begin{theorem}\label{t34}
{\rm \cite{GV}}
The multiplication $*$ and the bracket
$[.,.]$ induce a Gerstenhaber (or what is the same 1-Poisson)
algebra structure on the homology of the Hochschild complex $(\O,
\partial)$.~$\square$
\end{theorem}

\noindent{\bf Proof of Theorem~\ref{t34}:} The proof is deduced from the following
homotopy formulas.
$$
x*y-(-1)^{deg(x)deg(y)}y*x=(-1)^{deg(x)}(\partial(x\circ y)-
\partial x\circ y-(-1)^{deg(x)-1}x\circ \partial y).
\eqno(\numb)\label{eq35}
$$
The above formula proves the graded commutativity
of the multiplication $*$.

\vspace{2mm}
\noindent
$[x,y*z]-[x,y]*z-(-1)^{(deg(x)-1)deg(y)}y*[x,z]=$
$$
=(-1)^{deg(x)+deg(y)}
\bigl(\partial(x\{y,z\})-(\partial x)\{y,z\}-(-1)^{|x|}x\{\partial y,
z\}-(-1)^{|x|+|y|}x\{y,\partial z\}).
\eqno(\numb)\label{eq36}
$$
This formula proves the compatibility of the
bracket with the multiplication.~$\square$

\section{The little cubes operad}\label{s4}

Analogously to linear operads one can define topological operads,
{\it i.~e.}, collections $\{\O(n),n\ge 0\}$ of topological sets with

(i) an $S_n$-action on each $\O(n)$;

(ii) compositions
$$
\gamma:\O(l)\times\bigl( \O(m_1)\times\dots\times\O(m_l)\bigr)\to
\O(m_1+\dots +m_l);
$$

(iii) a unit element $id\in\O(1)$.

We assume the same associativity and symmetric group equivariance
requirements.

Evidently, the homology spaces $\{H_*(\O(n),\kk),n\ge 0\}$ 
over any field
$\kk$ form a graded $\kk$-linear operad. The same is true for any
commutative ring of coefficients, if the groups $H_*(\O(n),\Z)$,
$n\ge 0$, have no torsion.

Usually, when one considers the homology of a topological operad, one
inverses the grading, supposing that the $i$-th homology group
$H_i(\O(n),\kk)$ has the degree $-i$.

\vspace{1.5mm}

Historically one the first examples of topological 
operads are the little cubes
operads {$\LC_d=\{\LC_d(n),n\ge 0\}$},  $d\ge 1$, 
see~\cite{BV1},~\cite[Chapter~2]{BV2},~\cite[Chapter~4]{May}.
Here $\LC_d(n)$ denotes
the configuration space of $n$ disjoint cubes labelled by
$\{1,\dots,n\}$ in a unit cube. It is supposed that the faces of the
cubes are parallel to the corresponding faces of the unit cube.
The group $S_n$ acts by permutations of the cubes in configurations.
The element $id\in\LC_d(1)$ is the configuration of one cube coinciding
with the unit cube. The composition operations (ii) are insertions
of $l$ confi\-gu\-rations respectively of $n_1,\dots,n_l$ cubes into the
corresponding $l$ cubes of a configuration of $l$ cubes.

\begin{theorem}\label{t41}
{\rm \cite{Co}}
The homology groups of the operad $\LC_d$,
$d\ge 1$,
have no torsion and form the following graded linear
operad

1) if $d=1$, $\{H_{-*}(\LC_1(n),\kk),n\ge 0\}$ is isomorpic to the
associative algebras operad $\Assoc$,

2) if $d\ge 2$, $\{H_{-*}(\LC_d(n),\kk),n\ge 0\}$ is isomorphic to the
$(d-1)$-Posson algebras operad $\Poiss_{(d-1)}$.~$\square$
\end{theorem}

\begin{remark}\label{r42}
{\rm
The operad of associative algebras has a natural
filtration compatible with the operad structure. The
graded quotient associated with this filtration is the Poisson algebras
operad $\Poiss$.~$\square$
}
\end{remark}

We will not prove Theorem~\ref{t41}, however we will make some explanations.

First of all, note that the space $\LC_d(1)$ is contractible for any
$d\ge 1$. The homology classe of one point in this space
corresponds to $id=x_1\in\Poiss_{(d-1)}(1)$. The space $\LC_d(2)$ is
homotopy equivalent to the $(d-1)$-dimensional sphere $S^{d-1}$.
For $d\ge 2$ this gives us one operation of degree zero and one
operation of degree $1-d$. The first operation corresponds to the
multiplication
$x_1\cdot x_2\in\Poiss_{(d-1)}(2)$, the second one to the bracket
$[x_1,x_2]\in\Poiss_{(d-1)}(2)$.

Obviously, the space $\LC_d(n)$ is homotopy equivalent to
the configuration space of  collections of 
$n$ distinct points in $\R^d$, {\it i.~e.}, to
the space of injective maps $\{1,\dots,n\}\hookrightarrow\R^d$.
The latter space is an open everywhere dense subset of the vector space
$\R^{nd}$ of all maps $\{1,\dots,n\}\to \R^d.$ The complement
(called the {\it discriminant}) $\Delta_d(n)\subset\R^{nd}$ is a union
of $\frac{n(n-1)}{2}$ vector subspaces of codimension $d$.
Each of these subspaces $L_{i,j}$ corresponds to a pair of distinct
points $i,j\in\{1,\dots,n\}$ and consists of maps $\psi:\{1,\dots,n\}
\to\R^d$, such that $\psi(i)=\psi(j)$.

The non-complete partitions of the set $\{1,\dots,n\}$ (by the
{\it complete} partition we mean the partition into $n$ singletons)
are in one-to-one correspondence with the strata of the arrangement
$\Delta_d(n)=\bigcup_{i<j}L_{i,j}$. To a partition we assign a vector
subspace consisting of maps $\{1,\dots,n\}\to\R^d$ that glue the
points of each set in the partition. Following the ge\-ne\-ral theory
of arrangements the reduced homology groups $\tilde H_*(\R^{nd}\backslash
\Delta_d(n),\kk)$ can be decomposed into a direct sum, each summand
being assigned to some stratum of the arrangement,
see~\cite{GoMc, V2, V4, ZZ}.
In the
case $d\ge 2$ this decomposition is canonical (in the case $d=1$ it
depends on the choice of a component of $\R^{n}\backslash\Delta_1(n)).$
Let us assign to the
complete partition the degree zero homology group
$H_0(\R^{nd}\backslash\Delta_d(n)
,\kk)\simeq \kk$.
\medskip

\begin{proposition}\label{p43}
The above decomposition of
the homology groups $H_{-*}(\R^{nd}\backslash\Delta_d(n),\kk)$, $d\ge 2$, coincides
with the decomposition~(\ref{eq113}).~$\square$
\end{proposition}

\section{The first term of the Vassiliev auxiliary spectral sequence}\label{s5}

In the same way as the homology groups of the space of injective maps
$\{1,\dots,n\}\hookrightarrow\R^d$ are decomposed into a direct sum
by strata of the discriminant $\Delta_d(n)$ (see Section~\ref{s4} ---
Proposition~\ref{p43}), the first term of the Vassiliev auxiliary
spectral sequence is naturally decomposed into a direct sum in
which the summands are numbered by
equivalence classes of so called {\it $(A,b)$-configurations} defined below.

Let $A$ be a non-ordered finite collection of natural numbers
$A=(a_1,\dots,a_{\#A})$, any of which is not less than 2, and let $b$ be a
non-negative integer. Set $|A|:=a_1+\dots+a_{\#A}$. An {\it
$(A,b)$-configuration} is a collection of $|A|$ distinct points in
$\R^1$ separated into $\#A$ groups of cardinalities $a_1,\dots,
a_{\#A}$, plus a collection of $b$ distinct points in $\R^1$ (some
of which can coincide with the above $|A|$ points). For short
$(A,0)$-configuration are called simply {\it $A$-configurations}.
A map $\phi:\R^1\to\R^d$ {\it respects} an $(A,b)$-configuration, if
it glues together all points inside any of its groups of cardinalities
$a_1,\dots,a_{\#A}$, and the derivative $\phi'$ is equal to 0 at all
the $b$ last points of this configuration. For any $(A,b)$-configuration
the set of maps respecting it is an affine subspace in $\K$ of
codimension $d(|A|-\#A+b)$; the number $|A|-\#A+b$ is called the
{\it complexity} of the configuration. Two $(A,b)$-configurations are
called {\it equivalent} if they can be transformed into one another
by an orientation-preserving homeomorphism $\R^1\to\R^1$.

Consider any $(A,b)$-configuration $J$ of complexity $i$ and with
$j$ geometrically distinct points in $\R^1$. The stratum
consisting of all mappings $\R^1\to\R^d$  respecting at least one
$(A,b)$-configuration $J'$ equivalent to $J$, can be parameterized
by the space $S(J)$ of affine fiber bundle, whose base space is the
space $E^j$ of $(A,b)$-configurations $J'$ equivalent to $J$, and
the fiber over $J'$ is the affine space $\R^{\omega -di}$ of maps
respecting $J'$. Note that $E^j$ is contractible being an open cell
of dimension $j$. Therefore this fiber bundle can be trivialized:
$$
S(J)\simeq E^j\times\R^{\omega -di}.
\eqno(\numb)\label{eq51}
$$

\begin{remark}\label{r52}
{\rm
The corresponding stratum is not homeomorphic
to $S(J)$: one map $\R^1\to\R^d$  may respect two different
$(A,b)$-configurations equivalent to $J$, and therefore the stratum
has self-intersections.~$\square$
}
\end{remark}

The auxiliary spectral sequence computing the homology groups of the term
$\sigma_i\backslash\sigma_{i-1}$ in our filtration~(\ref{eq02})
uses those and only those $(A,b)$-configurations, that are of
complexity $i$. The geometrical meaning of the first differential
in the auxiliary spectral sequence is in how
the strata (or rather
the base spaces of the corresponding affine bundles~(\ref{eq51})) coresponding
to equivalence classes of $(A,b)$-configurations of complexity
$i$ bound to each other.

Let first  $d$ be even. For any equivalence class of
$(A,b)$-configurations of complexity $i$ and with $j$ geometrically
distinct points we will assign a star-partition $(\bar A,S)$
of the set
$\{1,\dots,j\}$ (see Definition~\ref{d119}) and therefore a subspace
$\BV_{(d-1)}(\bar A,S,j)$ of the operad $\BV_{(d-1)}$ (see
Decomposition~\ref{eq118}), where the partition $\bar A$ of the set
$\{1,\dots,j\}$  and the subset $S\subset\{1,\dots,j\}$ are
defined below.

\begin{definition}\label{d53}
{\rm
A {\it minimal component} of an
$(A,b)$-configuration is either one of its $b$ points, which
does not coincide with none of the $|A|$ points, or one of the
$\#A$ groups of points with all the stars contained there.~$\square$
}\end{definition}

For instance, the $(A,b)$-configuration on the Figure~\ref{f53'} has 3
minimal components consisting respectively of the following groups of points: 
\quad 1) $t_1$, $t_3$;\quad 
2) $t_2$, $t_4^*$, $t_5$;\quad
3) $t_6^*$.

\hspace*{32mm}
\PSbox{mcomp.pstex}{-10mm}{38mm}
\begin{picture}(10,10)
\put(172,20){\Large $*$}
\put(243,20){\Large $*$}
\end{picture}

\vspace{3mm}\vspace{1mm}

\parbox[b]{162mm}{
\centerline{(Figure \numb)\label{f53'}}
}
\vspace{3mm}

To any $(A,b)$-configuration with $j$ geometrically distinct
points we can assign a star-partition of the set $\{1,\dots,j\}$.
Let $t_1<t_2<\dots <t_j$ be the points in $\R^1$
of our $(A,b)$-configuration
($b$ of them are marked by stars).
The set of minimal components of this $(A,b)$-configuration defines
a partition $\bar A$ of the set $\{1,\dots,j\}$. We also have
the subset $S\subset\{1,\dots,j\}$ (of cardinality $b$) of indices,
corresponding to the points marked by stars.

Note that different equivalence classes of $(A,b)$-configurations
correspond to different star-partitions. But the correspondence
is far from being bijective. The star-partitions not corresponding
to $(A,b)$-configurations are those and only those, which contain
singletons not marked by a star.

\begin{theorem}\label{t54}
{\rm see \cite{T2}.}
The first term of the Vassiliev
auxiliary spectral sequence for even $d$ together with its first
differenrial is isomorphic to the subcomplex of the Hochschild complex
$(\BV_{(d-1)},\partial)$, linearly spanned by the summands of the
decomposition~(\ref{eq118}), corresponding to star-partitions, which don't
contain singletons not marked by a star. The grading corresponding
to the homology of the knot space $\K\backslash\Sigma$ is
minus the grading ``$deg$'' defined by~(\ref{eq33}).~$\square$
\end{theorem}

In particular, the theorem claims that the subspace in
$\BV_{(d-1)}$ spanned by the summands in question is invariant with
respect to the differential $\partial$.

Consider the decomposition of the space $\BV_{(d-1)}=\bigoplus_
{i\ge 0}E_i$, where $E_i$ is the sum over all star-partitions
having exactly $i$ singletons not marked by a star. The filtration
$F_0\supset F_1\supset F_2\dots$, with $F_i:=\bigoplus_{j\ge i} E_j$,
is compatible with the differential $\partial$. Note, that the complex
$(E_0,\partial)$ is exactly the complex from Theorem~\ref{t54}.

\begin{proposition}\label{p55}
{\rm \cite{T2}}
The Hochschild complex $(\BV_{(d-1)},
\partial)$ ($d$ is even) is a direct sum of the complexes $(E_0,
\partial)$ and $(F_1,\partial)$. The first complex $(E_0,\partial)$
is homology (and even homotopy) equivalent to $(\BV_{(d-1)},
\partial)$ ; the second one $(F_1,\partial)$ is acyclic (and even
cont\-ractible).~$\square$
\end{proposition}

An analogous statement holds for the complexes $(\Poiss_{(d-1)},\partial)$,
$d$
being any integer number. The subcomplexes spanned
by the summands of the decomosition~(\ref{eq113}) cor\-res\-ponding
to partitions non-containing singletons will be called the {\it normalized
Hochschild complexes} and denoted by
$(\Poiss_{(d-1)}^{Norm},\partial)$. The {\it normalized Hochschild
complexes} $(E_0,\partial)$ for the operads $\BV_{(d-1)}$,
$d$ being even, will be denoted also 
by $(\BV_{(d-1)}^{Norm},\partial)$.

\begin{corollary}\label{cor56}
{\rm
The first term of the main spectral sequence
(for even $d$ and for any commutative ring $\kk$ of coefficients)
is isomorphic to the Hochschild homology (with inversed grading)
 of the Batalin-Vilkovisky
operad $\BV_{(d-1)}$.
}~$\square$
\end{corollary}

\begin{remark}\label{r57}
{\rm
(due to M.~Kontsevich) The operator $\delta\in\BV
_{(d-1)}(1)$ has a natural geometrical interpretation as  the Euler
class $\delta_{d-1}^E\in H_{d-1}(SO(d))$ of the special orthogonal
group $SO(d)$, see Section~\ref{s7}.~$\square$
}
\end{remark}

If $d$ is odd, then the first term of the auxiliary spectral sequence is
very similar to the normalized Hochschild complex $(\BV^{Norm},\partial
)$, but it does not correspond to any operad. A description of the
obtained complex, which I called the {\it complex of bracket
star-diagrams}, see in Section~\ref{s8}. 

In~\cite{T2}  the following theorem is  proved. For even $d$ this is an
immediate corollary of Theorem~\ref{t54}.

\begin{theorem}\label{t58}
For any $d\ge 3$ the subspace of the first
term of the auxiliary spectral sequence linearly spanned only by
the summands corresponding to $A$-configurations forms a subcomplex
isomorphic to the normalized Hochschild complex $(\Poiss_{(d-1)}^
{Norm},\partial)$ with inversed grading.~$\square$
\end{theorem}

There is a very nice way to simplify the computations of the second
term of the auxiliary spectral sequence. This construction works
both for $d$ even and odd. Consider the nor\-ma\-li\-zed Hochschild
complex $(\Poiss_{(d-1)}^{Norm},\partial)$ and consider the quotient of the
space of
this complex by  the ``neighboring commutativity relations'' ---
in  other words, for any $i=1,2,  \dots$ we set $x_i$ and $x_{i+1}$
to commute. For example the element $[x_1,x_3]\cdot[x_2,[x_4,x_5]]\in
\Poiss_{(d-1)}(5)$ is equal to zero modulo these relations, because
it contains $[x_4,x_5]$. The space of relations
is invariant with respect to the differential, thus the quotient space
has the structure of a quotient complex. Denote this quotient
complex  by $(\Poiss_{(d-1)}^{zero},\partial)$.

\begin{theorem}\label{t59}
{\rm \cite{T2,T3}}
For any commutative ring $\kk$
of coefficients the space $\Poiss_{(d-1)}^{zero}$ is a free
$\kk$-module. The homology space of the complex $(\Poiss_{(d-1)}^{zero},
\partial)$ is isomorphic to the first term of the main spectral
sequence (=to the second term of the auxiliary spectral
sequence).~$\square$
\end{theorem}

\noindent{\bf Idea of the proof:} Consider the filtration in the 
first term (of the auxiliary spectral sequence)
 by the number of minimal components of corresponding
$(A,b)$-configurations.
The associated spectral sequence degenerates
in the second term, because its first term is concentrated on the
only line corresponding to $A$-configurations.~$\square$

\vspace{2.5mm}

Before describing  the relation of the Hochschild homology space
 of
$\Poiss_{(d-1)}$ defined over $\Q$ with
the first term of the main spectral sequence also defined over $\Q$
(see Theorem~\ref{t510}) we will give some explanations.

In Section~\ref{s8} a structure of  a differential graded
cocommutative bialgebra on the first term of the auxiliary spectral
sequence is defined. According to Theo\-rems~\ref{t54}
and~\ref{t58} such a structure is defined also
on the normalized Hochschild complexes $(\Poiss^{Norm},
\partial)$, $(\Gerst^{Norm},\partial)$, $(\BV^{Norm},\partial)$. The
neighboring commutativity relations respect multiplication and
comultiplication (but not the bracket~\eqref{bracket}), so the complexes
$(\Poiss_{(d-1)}^{zero},\partial)$ are differential bialgebras.
A motivation of the existence of such a  structure is that the
space of long knots is an $H$-space (has a homotopy associative
multiplication), see~\cite{T1,T2,T3};  therefore its homology
and cohomology spaces
over any field $\kk$ are mutually dual graded respectively
cocommutative and commutative bi\-al\-geb\-ras. For $\kk=\Q$ these
(co)homology bialgebras are bicommutative, see~\cite{T2}, \cite{T3}. This fact
supports 
Conjecture~\ref{con03}, because applying Theorem~\ref{t34} to 
the operads of Poisson, Gerstenhaber or Batalin-Vilkovisky algebras
over any field $\kk$, we obtain that their  Hochschild homology
bialgebras are bicommutative. It follows from the Milnor
theorem, see~\cite{MM}, that for $\kk=\Q$
 these bialgebras are polynomial. The space
of pri\-mi\-tive elements is the space of their generators. In fact,
the bracket~\ref{bracket} for these operads preserves the spaces of primitive
elements.
We conclude that as Gerstenhaber algebras the Hochschild homology 
spaces
of these operads are symmetric algebras of the Lie algebras of
primitive elements (see Example~\ref{ex110} with $d=1$).

\begin{theorem}\label{t510}
{\rm \cite{T2,T3}}
As a graded
bialgebra the first term of the main
spectral
sequence over $\Q$ is isomorphic to the Hochschild homology bialgebra
(with the inversed grading) of
the $\Q$-linear operad $\Poiss_d$  factorized 

1) for even $d$: by one odd primitive generator $[x_1,x_2]$ (of degree $3-d$);

2) for odd $d$: by one even primitive generator $[x_1,x_2]$ (of degree $3-d$)
and one odd primitive generator $[[x_1,x_3],x_2]$ (of degree $5-2d$).~$\square$
\end{theorem}

\section{The bialgebra of chord diagrams}\label{s6}

The {\it bialgebra of chord diagrams},
the dual to the associated
quotient bialgebra of Vassiliev knot
invariants, was intensively
studied in the last decade. In this section we give an
interpretation of the bialgebra of chord diagrams as a part of the
Hochschild homology algebra  of the Poisson algebras operad.

Consider the normalized Hochschild complex $(\Poiss_{(d-1)}^{Norm},\partial)$.
In this complex one can define a bigrading by the complexity
$i$ and by the number $j=|A|$ of geometrically distinct points of
the corresponding $A$-configurations. The differential $\partial$ is
of bidegree (0,1). If the first grading $i$ is fixed, then 
the number $j$
varies from $i+1$ to $2i$. So any element of the bigrading $(i,2i)$ belongs
to the kernel of $\partial$. The case $j=2i$ corresponds to the minimal
possible dimension    of non-trivial homology classes 
of the space of long knots for the complexity $i$
fixed. For example, if $d=3$, then this dimension is equal to
$(d-1)i-j=2i-j=0$.
The part of the Hochschild homology groups,
that lies in the bigradings
$(i,2i)$, $i\ge 0$, will be called the {\it bialgebra of
chord diagrams} if
$d$ is odd, and the {\it bialgebra of chord superdiagrams}
if $d$ is even.
Any product of brackets in $\Poiss_{(d-1)}^{Norm}(2i)
\subset Poiss_{(d-1)}(x_1,
\dots,x_{2i})=
S^* Lie_{(d-1)}(x_1,\dots,x_{2i})$ of 
bidegree $(i,2i)$ is the  product
of $i$ brackets, each of wich contains exactly 2 generators. Thus,
any such product of brackets can be depicted as $2i$ points on the line
$\R^1$, that are decomposed into $i$ pairs and connected by a chord
inside each pair. For example, $[x_3,x_5]\cdot[x_4,x_1]\cdot [x_2,x_6]$
is assigned to the diagram

\hspace*{45mm}
\PSbox{chd.pstex}{15mm}{30mm}
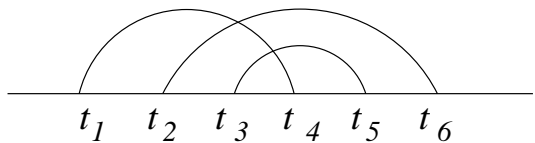
\begin{picture}(10,10)
\end{picture}

\vspace{3mm}\vspace{1mm}

\parbox[b]{165mm}{
\centerline{(Figure \numb)\label{f61}}
}
\vspace{3mm}

The so called {\it  4-term relations} arise as the differential
$\partial$ of products of brackets, in which all brackets except one
are ``chords'', and the only non-chord is a bracket on three elements. In
other words, these products of brackets correspond to $A$-configurations,
with $A=(3,2,2,\dots,2)$.

\begin{remark}\label{r62}
{\rm
Sometimes one considers the bialgebra of chord
diagrams factorized not only by 4-term relations, but also by 1-term
relations. The latter relations arise (when we consider the whole first
term of the auxiliary spectral sequence)
as the differential of diagrams corresponding to $(A,b)$-configurations,
with $A=(2,\dots,2)$, $b=1$, and the only star  does not  coincide with
none of $|A|$ points.~$\square$
}
\end{remark}

In the case of odd $d$ we need to take into account orientations of
chords (since $[x_{i_1},x_{i_2}]=-[x_{i_2},x_{i_1}]$). This definition
does not coincide with the standard one, where these orientations
are not important. This discordance  of  definitions can be
easily eliminated. Consider the Hochschild complex
$(\Poiss_{(d-1)},\partial)=\bigl(\bigoplus_{n\ge 0}\Poiss_{(d-1)}(n),
\partial\bigr)$, and replace each space $\Poiss_{(d-1)}(n)$ by its
tensor product with  the  one-dimensional
sign representation
$sign$ of the symmetric group $S_n$. This can be done in the following way.  We
take a free $(d-1)$-Poisson algebra $Poiss_{(d-1)}
(x'_1,\dots,x'_n)=S^*Lie(x'_1,\dots,x'_n)$
with ge\-ne\-ra\-tors $x'_1,\dots,x'_n$ of degree one instead of the analogous
algebra $Poiss_{(d-1)}(x_1,\dots,x_n)\supset\Poiss_{(d-1)}(n)$ with
generators $x_1,\dots,x_n$ of degree zero. Afterwards we consider
the  subspace $\Poiss'_{(d-1)}(n)\subset Poiss_{(d-1)}(x'_1,\dots,x'_n)$
spanned by products of brackets con\-tai\-ning each generator $x'_i$
exactly once. Obviously, the $S_n$-module $\Poiss'_{(d-1)}(n)$ is
isomorphic to $\Poiss_{(d-1)}(n)\otimes sign$.
Defining properly the differential, see Section~\ref{s8}, we obtain 
another version
of the
Hochschild complex for the
operads $\Poiss_{(d-1)}$.
This new version is interpreted in the geometry of the discriminant
$\Sigma$ as introducing an orientation of the spaces $E^j$ of the strata~(\ref{eq51}) not according to the usual order
$t_1<t_2<\dots<t_j$ of the points on the line $\R^1$,
but according to the order, that was in our product of brackets. 
For
instance, for the element $[x_3,x_4]\cdot [x_5,[x_2,x_1]]$ the orientation of
the corresponding space $E^5=\{t_1<t_2<t_3<t_4<t_5\}$ of 
$(3,2)$-configurations equivalent to
configuration~\eqref{f63} is according to the order $(t_3,t_4,t_5,t_2,t_1)$.

\hspace*{50mm}
\PSbox{aconf.pstex}{15mm}{38mm}
\begin{picture}(10,10)
\end{picture}

\vspace{3mm}\vspace{1mm}

\parbox[b]{165mm}{
\centerline{(Figure \numb)\label{f63}}
}
\vspace{3mm}

An advantage of the new Hochschild complexes (for
operads $\Poiss_{(d-1)}$) is a simpler rule of signs in the definition
of the differential, see Section~\ref{s8}.

\vspace{1.5mm}

Let $\O$ be a graded linear operad, equipped with a morphism from
the operad $\Assoc$, then any $\O$-algebra $A$ is an associative algebra
because of the following morphisms:
$$
\Assoc\to\O\to \Endom(A).
$$
On the other hand, the map $\O\to \Endom(A)$ defines a morphism of
Hochschild complexes:
$$
(\O,\partial)\to(\Endom(A),\partial).
$$
Therefore the classes in the Hochschild homology of $\O$ can be
considered as characteristic classes of the Hochschild
cohomology of $\O$-algebras (considered as associative algebras).
An interesting question is whether all the classes in the Hochschild
homology of $\BV$, $\Poiss$, $\Gerst$, $\BV_{(d-1)}$, $\Poiss_{(d-1)}$
have a non-trivial realization as characteristic classes.

Consider Example~\ref{ex110}, where  we take $S^*(\g[d-1])$
as a $(d-1)$-Poisson algebra.
The Hochschild homology space
of a polynomial algebra
is well known. It is the space of  polynomial  polyvector fields on the
space of generators. In our case the space of generators is $\g[d-1]$.
According to the grading rule, we get that the Hochschild homology
space
of $\Poiss_{(d-1)}$ of bigrading $(i,j)$ is mapped to the space
of homogenous  degree  $i$ $j$-polyvector fields, 
{\it i.~e.},
of  expressions of the form
$$
\sum_{{q_1,\dots,q_j}\atop{p_1,\dots,p_i}} A_{p_1\dots p_i}^{q_1\dots q_j}
x^{p_1}\dots x^{p_i}\frac\partial
{\partial x^{q_1}}\wedge\dots\wedge\frac\partial{\partial x^{q_j}}.
\eqno(\numb)\label{eq64}
$$
If $d$ is odd (resp. even), the tensor $A_{p_1\dots p_i}^{q_1\dots q_j}$
is symmetric (resp. antisymmetric) with respect to the indices $p_1,\ldots,p_i$
and antisymmetric (resp. symmetric) with respect to the indices 
$q_1,\ldots,q_j$.
\vspace{2.5mm}

{\bf Problem \numb.\label{p65}} Find explicitly
$A_{p_1\dots p_i}^{q_1\dots q_j}$
via the bracket on $\g$, for example, in the case $j=2i$ of chord
diagrams.~$\square$
 
\vspace{2.5mm}

The answer to Problem~\ref{p65}  can be related with
the invariant tensors, defined for Casimir Lie algebras by chord diagrams,
see~\cite{BN1, HV, Vai}.

\bigskip

Also a problem which looks interesting is to find the Gerstenhaber subalgebra
in the Hochschild homology of the operads $\Poiss$, $\Gerst$ generated
by the space of chord (super)diagrams. For instance, in the
case of the operad $\Gerst$ there is an element of bigrading (3,5),
that cannot be obtained from chord diagrams by means of
the multiplication and the bracket~\eqref{bracket}, see~\cite{T1, T2}.
In the case
of the operad $\Poiss$ we do not know such an example.

\section{Operads of turning balls}\label{s7}

In this section we 
clarify the difference between the cases of odd and even $d$
and  give a geometrical interpretation of the operad $\BV_{(d-1)}$
(for $d$  even). In particular we explain
Remark~\ref{r57} (due to 
M.~Kontsevich). Theorems~\ref{t71}, \ref{t72} given below
are classical, see~\cite{G},~\cite{B}.

\vspace{1.5mm}

For any $d\ge 2$ let us introduce  the  topological {\it operad
$\TB=\{\TB_d(n),n\ge 0\}$ of turning balls}. The space $\TB_d(n)$ is
put to be the configuration space of $n$ mappings of the  unit ball
$B^d\subset\R^d$ into itself. These $n$ mappings are supposed

1) to be injective, preserving the orientation and the ratio of
distances;

2) to have disjoint images.

\vspace{1mm}

The space $\TB_d(n)$ can be evidently identified with the 
direct product of the
$n$-th power $(SO(d))^n$ of the special orthogonal group with the configuration
space of $n$ disjoint balls in the unit ball $B^d$. The last space is homotopy
equivalent to the $n$-th component $\LC_d(n)$ of the little cubes operad. The
composition operations, the symmetric group actions and the identity element
are defined analogously to the case of the little cubes operad.

\begin{theorem}\label{t71}
{\rm \cite{G}}
The homology $\{H_{-*}(\TB_2(n),\Z),n\ge 0\}$
of the turning discs operad (balls of dimension 2) is the
Batalin-Vilkovisky algebras operad $\BV$.~$\square$
\end{theorem}

Let us describe the operad $\{H_{-*}(\TB_d(n),\Q),n\ge 0\}$ for any
$d\ge 2$. Note, that the space of  unary operations is the
homology algebra $H_{-*}(SO(d),\Q)$.

\begin{theorem}\label{t72}
{\rm \cite{B}}
The homology bialgebra $H_*(SO(d),\Q)$ is the exterior
algebra on the following pri\-mi\-tive generators:

1) case $d=2k+1$: generators $\delta_3,\delta_7,\dots,\delta_{4k-1}$
of degree $3,7,\dots,4k-1$ respectively;

2) case $d=2k$: generators $\delta_3,\delta_7,\dots,\delta_{4k-5},
\delta_{2k-1}^E$ of degree $3,7,\dots,4k-5$ and $2k-1$
respectively.~$\square$
\end{theorem}

The generators $\delta_{4i-1}\in H_{4i-1}(SO(d),\Q)$, $d\ge 2i+1$, are
called the {\it Pontriagin classes}, the generator $\delta_{2k-1}^E\in
H_{2k-1}(SO(2k),\Q)$  is called the {\it Euler class}.
The Pontriagin classes (and the subalgebra generated by them) lie
in the kernel of the map in homology $H_*(SO(d),\Q)\to H_*(S^{d-1},\Q)$
induced by the natural
projection
$$
SO(d)\stackrel{SO(d-1)}\longrightarrow S^{d-1}.
$$
The Euler class is sent to the canonical class of dimension $d-1$
in the homology of the sphere $S^{d-1}$.

Now we are ready to describe the operads
$\{H_{-*}(\TB_d(n),\Q),n\ge 0\}$, $d\ge 2$.
We will say which objects are algebras over these operads.

\begin{theorem}\label{t73}
{\rm \cite{T2}}
Algebras over the operad
$\{H_{-*}(\TB_d(n),\Q),n\ge 0\}$ are

1) for even $d$: $(d-1)$-Batalin-Vilkovisky algebras, where the
operator $\delta$ is $\delta_{d-1}^E$;

2) for odd $d$: $(d-1)$-Poisson algebras.

Furthemore these algebras are supposed to have $[\frac{d-1}2]$
(where ``$[\,.\,]$'' denotes the integral part) mutually super-comuting
differentials $\delta_3,\delta_7,\dots,\delta_{4[\frac{d-1}2]-1}$
(of degree $-3,-7,\dots,1-4[\frac{d-1}2]$) of the Batalin-Vilkovisky
(resp. Poisson) algebra structure. This means, that for any elements
$a,b$ of the algebra and $1\le i\le [\frac{d-1}2]$, one has

(i) $\delta_{4i-1}([a,b])=[\delta_{4i-1}a,b]+(-1)^{\tilde a +d-1}[a,
\delta_{4i-1}b]$;

(ii) $\delta_{4i-1}(a\cdot b)=(\delta_{4i-1} a)\cdot b +(-1)^{\tilde{a}}
a\cdot(\delta_{4i-1} b)$;

(iii) only for even $d$:
\quad $\delta_{4i-1}\delta_{d-1}^Ea=-\delta_{d-1}^E\delta_{4i-1}a.
$~$\square$
\end{theorem}

\section{Complexes of bracket star-diagrams}\label{s8}

In this section we describe the first term of the auxiliary spectral sequence
together with its first differential. The corresponding complex will be called
{\it complex of bracket star-diagrams}. In the case of even $d$ this complex
is isomorphic to the normolized Hochschild complex of the Batalin-Vilkovisky 
algebras operad $(\BV^{Norm},\partial)$.

\subsection{Case of odd $d$}\label{s81}

Let us fix an $(A,b)$-configuration $J$. Let 
$t_\alpha,$ $\alpha\in${\LARGE $\alpha$} (resp. $t^*_\beta$, 
$\beta\in${\LARGE $\beta$}) be the points of $J$ on the line $\R^1$ that do
not have stars (resp. that do have stars). Consider the free Lie super-algebra
with the even bracket and with odd generators $x_{t_\alpha}$, 
$\alpha\in${\LARGE $\alpha$}, $x_{t^*_\beta}$, $\beta\in${\LARGE $\beta$}.
We will take the symmetric  algebra (in the super-sens) of the space of this
Lie super-algebra. In the obtained space we will consider the subspace
$\BSD(J)$ spanned by the (products of brackets), where each minimal component
of $J$ is presented by one bracket, containing only generators indexed
by the points of this minimal component and containing each such generator 
exactly once.

Such products of brackets will be called {\it bracket star-diagrams}.

\begin{example}\label{ex81}
{\rm
The space $\BSD(J)$ of the bracket star-diagrams corresponding to the 
$(A,b)$-configuration $J$ of the Figure~\ref{f53'} is two-dimensional.
The diagrams
$$[x_{t_1},x_{t_3}]\cdot[[x_{t_2},x_{t_4^*}],x_{t_5}]\cdot x_{t_6^*},\,\,
[x_{t_1},x_{t_3}]\cdot[[x_{t_2},x_{t_5}]x_{t_4^*}]\cdot x_{t_6^*}$$
form a basis in this space.~$\square$
}
\end{example}

If two bracket star-diagrams can be transformed into one another by an 
orientation preserving homeomorphism $\R^1\to\R^1$, then they are set 
to be equal. For any equivalence class ${\bf J}$ of $(A,b)$-configurations 
we define the space $\BSD({\bf J})$ as the space $\BSD(J)$, where $J$ is any
element of ${\bf J}$.

The {\it space of bracket star-diagrams} is defined as the direct sum of the
spaces $\BSD({\bf J})$ over all equivalence classes ${\bf J}$ of 
$(A,b)$-configurations.

The complexity $i$ and the number $j$ of geometrically distinct points
of the corresponding $(A,b)$-configurations define the bigrading $(i,j)$
on the space of bracket star-diagrams. Remind that the complex of bracket 
star-diagrams (both for $d$ odd and even) is supposed to compute the first
term $E^1_{p,q}$ of the main spectal sequence, whose $(p,q)$ coordinates are 
expressed as follows:
$$p=-i,
$$
$$q=id-j.
$$
The corresponding homology degree
of the space of knots $\K\backslash\Sigma$ is 
$$p+q=i(d-1)-j.
\eqno(\numb\label{eq82})
$$

Note that the first term of the main spectral sequence is non-trivial only
in the second quadrant $p\le 0$. The inequality $j\le 2i$ for 
$(A,b)$-configurations provides the condition $q\ge(2-d)p$.

\hspace*{50mm}
\PSbox{spposl.pstex}{75mm}{80mm}
\begin{picture}(10,10)
\put(-230,100){$\tan\alpha=d-2$}
%\put(243,20){\Large $*$}
\end{picture}

\nopagebreak

\vspace{3mm}\vspace{1mm}

\nopagebreak

\parbox[b]{160mm}{
\centerline{(Figure \numb)\label{f82'}}
}
\vspace{3mm}

The bialgebra of chord (super)diagrams occupies the diagonal $q=(2-d)p$.

To describe the differential $\partial$ on the space of bracket star-diagrams
we will need some complementary definitions and notations.

\begin{definition}\label{d83}
{\rm
Let us permit to $(A,b)$-configurations, with $A=(a_1,\dots,a_{\#A})$, to 
have $a_i=1$; we demand also that one-element sets should never coincide with 
stars. These $(A,b)$-configurations will be called {\it generalized
$(A,b)$-configurations}. The generalized $(A,b)$-configurations that
are not $(A,b)$-configurations in the usual sens (that have $a_i=1$)
will be called 
{\it special} generalized $(A,b)$-configurations.~$\square$
}
\end{definition}

Analogously we define the {\it space of (special) generalized 
bracket  star-diagrams}.

\begin{example}\label{ex84}
$[x_{t_1},x_{t_3}]\cdot x_{t_2}$ is a special generalized star-diagram.~$\square$
\end{example}

Note that the space of generalized bracket star-diagrams is the direct sum 
of the space of bracket star-diagrams with the space of special generalized
bracket star-diagrams.

\begin{definition}\label{d85} 
{\rm We say that a (generalized) 
$(A,b)$-configuration $J$ can be {\it inserted} in a point $t_0^{(*)}$ of 
another 
(generalized) $(A',b')$-configuration $J'$, if $J$ does not have common points with $J'$ except possibly the point $t_0^{(*)}$.~$\square$
}
\end{definition}

\begin{definition}\label{d86}
{\rm
We say that a (generalized) bracket star-diagram can be {\it inserted} 
in the point 
$t_0^{(*)}$  of another (generalized) bracket star-diagram, if it is the
case for their $(A,b)$-configurations.~$\square$
}
\end{definition}

Let $A$ and $B$ be two  generalized bracket star-diagrams, such that $A$
can be inserted in the point $t_0$ (or $t_0^*$) of $B$, define the element
$B|_{x_{t_0}=A}$ (resp. $B|_{x_{t^*_0}=A}$) of the space  of 
generalized bracket star-diagrams. Up to a sign $B|_{x_{t^{(*)}_0}=A}$
is defined by replacing $x_{t^{(*)}_0}$ (in the diagram $B$) for $A$.
The sign is defined as $(-1)^{(\tilde{A}-1)\times n}$, where $\tilde A$ is the
parity of $A$ (the parity of the number of geometrically distinct points),
$n$ is the number of generators of the form $x_{t_\alpha },x_{t^*_\beta }$
before $x_{t^{(*)}_0}$ in $B$. In other words : we put the bracket containing
$x_{t_0^{(*)}}$ on the first place, then by means of antisymmetry relations
we put $x_{t_0^{(*)}}$ on the first place in the bracket (and therefore in 
the diagram); we replace $x_{t^{(*)}_0}$ for $A$; and we do all these
manipulations in the inverse order. It is easy to see that these two 
definitions give the same sign.

\begin{example}\label{ex87}
$$[x_{t_2}x_{t_3^*}]\cdot [x_{t_1^*}x_{t_0}]\bigl| _{x_{t_0}=[x_{t_4}x_{t_5}]
\cdot x_{t_6^*}}=(-1)^{(3-1)\cdot 3}[x_{t_2}x_{t_3^*}]\cdot
[x_{t_1^*},[x_{t_4}x_{t_5}]\cdot x_{t_6^*}]. ~\square$$
\end{example}

Note that if $A$ has more than 1 minimal components, then the element 
$B|_{x_{t^{(*)}_0}=A}$ contains multiplications inside brackets. Therefore 
it is no more a (generalized) bracket star-diagram. To express this element
as a sum of (generalized) bracket star-diagrams we will use the 
formula~\eqref{eq111}. 

Now we are ready to define the differential $\partial$ on the space of bracket
star-diagrams.

Let $A$ be a bracket star-diagram, and let $t_\alpha$ be one of its points
without a star, then we define
$$
\partial_{t_\alpha}A:=P\bigl(A|_{x_{t_\alpha }=x_{t_{\alpha -}}
\cdot x_{t_{\alpha +}}}\bigr),
\eqno(\numb\label{eq88})$$
where $P$ is the projection of the space of generalized bracket star-diagrams
on the space of bracket star-diagrams,
that sends the space of special generalized bracket star-diagrams to zero;
the points $t_{\alpha_-},t_{\alpha_+}\in\R^1$ are respectively $t_\alpha-
\epsilon$ and $t_\alpha+\epsilon$ for a very small $\epsilon>0$.

\begin{remark}\label{r89}
{\rm
The formula~\eqref{eq88} can be made more precise:
$$\partial_{t_\alpha }A+(x_{t_{\alpha -}}-x_{t_{\alpha +}})\cdot
A=A|_{x_{t_\alpha }=x_{t_{\alpha -}}\cdot x_{t_{\alpha +}}}.~\square
\eqno(\numb\label{eq810})$$
}
\end{remark}

Let now $t_\beta^*$ be one of the points of $A$ having a star. We define
$$
\partial_{t_\beta ^*}A:=P\bigl(
A|_{x_{t^*_\beta }=x_{t_{\beta -}}\cdot x_{t^*_{\beta
+}}+x_{t^*_{\beta -}}\cdot x_{t_{\beta +}}+[x_{t_{\beta -}},x_{t_{\beta
+}}]}\bigr),
\eqno(\numb\label{eq811})$$
where $P$ is the same projection, the points $t_{\beta-}^{(*)}$, 
$t_{\beta +}^{(*)}$
are respectively $t_\beta^*-\epsilon$ and $t_\beta^*+\epsilon$ for a very
small $\epsilon>0$.

\begin{remark}\label{r812}
{\rm
The formula~\eqref{eq810} can be made more precise:
$$
\partial_{t^*_\beta }A+(x_{t_{\beta -}}-x_{t_{\beta +}})\cdot 
A=A|_{x_{t^*_\beta
}=x_{t_{\beta -}}\cdot x_{t^*_{\beta +}}+x_{t^*_{\beta -}}\cdot x_{t_{\beta
+}}+ [x_{t_{\beta -}},x_{t_{\beta +}}]}.~\square
\eqno(\numb\label{eq813})$$
}
\end{remark}

The differential $\partial$ on the space of bracket star-diagrams is the 
sum of the operators $\partial_{t_\alpha}$ and $\partial_{t_\beta^*}$ over
all points $t_\alpha$, $\alpha\in${\LARGE $\alpha$}, and $t_\beta^*$,
$\beta\in${\LARGE $\beta$}, of the corresponding $(A,b)$-configurations:
$$
\partial =\sum_{\alpha\in
{\displaystyle \alpha }}\partial_{t_\alpha}+\sum_{\beta \in{\displaystyle
\beta }}\partial_{t_\beta^*}.
\eqno(\numb\label{eq814})$$

It is easy to see that $\partial^2=0$.

\begin{remark}\label{r815}
$$\partial A= \left( \sum_{\alpha\in
{\displaystyle \alpha }}
A|_{x_{t_\alpha }=x_{t_{\alpha -}}\cdot x_{t_{\alpha +}}}\right)
+ \left( \sum_{\beta \in{\displaystyle
\beta }}A|_{x_{t^*_\beta
}=x_{t_{\beta -}}\cdot x_{t^*_{\beta +}}+x_{t^*_{\beta -}}\cdot x_{t_{\beta
+}}+ [x_{t_{\beta -}},x_{t_{\beta +}}]}\right) -
(x_{t_-}-x_{t_+})\cdot A,
\eqno(\numb\label{eq816})
$$
{\rm
where $t_-$ (resp. $t_+$) is less (resp. greater) than all the points of 
the diagram $A$ on the line $\R^1$.~$\square$
}
\end{remark}

\subsection{Case of even $d$}\label{s82}

Let us fix an $(A,b)$-configuration $J$ and consider the free Lie super-algebra
with the even bracket and with the even generators $x_{t_\alpha}$,
$\alpha\in${\LARGE $\alpha$}, and the odd generators $x_{t_\beta^*}$,
$\beta\in${\LARGE $\beta$}, where $t_\alpha$, $\alpha\in${\LARGE
$\alpha$}, and $t_\beta^*$, $\beta\in${\LARGE $\beta$}, are the  points
of our $(A,b)$-configuration $J$. Let us take the exterior
algebra (in the super-sens) of the space of this Lie super-algebra. By
convention the parity of an element $A=A_1\wedge\ldots\wedge A_k$ is 
$\tilde A=\tilde A_1+\ldots+\tilde A_k +k-1$ the sum of the parities of 
$A_i$, $1\le i\le k$, plus the number $k-1$ of the exterior product signs.
In the obtained space we will consider the subspace $\BSD(J)$ spanned
by the analogous products of brackets (see the previous subsection~\ref{s81}).
These products of brackets will be also called {\it bracket star-diagrams}.
The {\it space of bracket star-diagrams} is defined analogously to the case
of odd $d$. We also accept Definitions~\ref{d83},~\ref{d85},~\ref{d86}.
Note that the parity of a (generalized) bracket star-diagram is the 
number of stars plus 
the number of the exterior product signs. This parity is opposite to
the parity of the corresponding homology degree~\eqref{eq82}.

Let $A$ and $B$ be two generalized bracket star-diagrams, such that $A$ 
can be inserted in the point $t_0^{(*)}$ of $B$. Let us define 
$B|_{x_{t_0^{(*)}}=A}$. To do this we replace  $x_{t^{(*)}_0}$ in $B$ by  $A$,
an we multiply the obtained expression by 
$(-1)^{(\tilde{A}-\epsilon_0)\times (n_1+n_2)}$, where $\epsilon_0$ is 
equal to zero (resp. to one) if the point $t_0$ has no star
(resp. if the point $t_0^*$ has a star); $n_1$ (resp. $n_2$)
is the number of the exterior product signs (resp. of the generators 
corresponding to stars) before $x_{t_0^{(*)}}$ in $B$.

\begin{example}\label{ex817}
$$[x_{t_2}x_{t_3^*}]\wedge [x_{t_1^*}x_{t_0}]\bigl| _{x_{t_0}=[x_{t_4}x_{t_5}]
\wedge x_{t_6^*}}=(-1)^{(2-0)\cdot (1+2)}[x_{t_2}x_{t_3^*}]\wedge
[x_{t_1^*},[x_{t_4}x_{t_5}]\wedge x_{t_6^*}]. ~\square$$
\end{example}

Let $A$ be a bracket star-diagram, let $t_\alpha$ be a point of $A$ 
without a star. Define
$$
\partial_{t_\alpha}A:=P\left(A|_{x_{t_\alpha}=x_{t_{\alpha-}}\wedge
x_{t_{\alpha+}}}\right),
\eqno(\numb\label{eq818})$$
where $P$ is the projection from the space of generalized bracket 
star-diagrams to the space of bracket star-diagrams.

\begin{remark}\label{r819}
{\rm
The formula~\eqref{eq818} can be made more precise:
$$
\partial_{t_\alpha}A+(x_{t_{\alpha-}}-x_{t_{\alpha+}})\wedge
A=A|_{x_{t_\alpha}=x_{t_{\alpha-}}\wedge x_{t_{\alpha+}}}.~ \square
\eqno(\numb\label{eq820})$$
}
\end{remark}

Let $t_\beta^*$ be a point of $A$ having a star, we define
$$
\partial_{t^*_\beta} A:=P\left(A|_{x_{t^*_\beta}=x_{t_{\beta-}}\wedge
x_{t_{\beta+}^*}-x_{t_{\beta-}^*} \wedge x_{t_{\beta+}}-[x_{t_{\beta-
}},x_{t_{\beta+}}]}\right).~\square
\eqno(\numb\label{eq821}) $$

\begin{remark}\label{r822}
{\rm
The formula~\eqref{eq821} can be made more precise:
$$
\partial_{t^*_\beta}A+(x_{t_{\beta-}}-x_{t_{\beta+}})\wedge
A=A|_{x_{t^*_\beta}=x_{t_{\beta-}}\wedge x_{t^*_{\beta+}}-
x_{t^*_{\beta-}}\wedge x_{t_{\beta+}}-[x_{t_{\beta-
}},x_{t_{\beta+}}]}.~\square
\eqno(\numb\label{eq823})$$
}
\end{remark}

The differential $\partial$ on the space of bracket star-diagrams is
defined by the formula~\eqref{eq814} analogously to the case of odd $d$.

\begin{remark}\label{r824}
$$\partial A= \left( \sum_{\alpha\in
{\displaystyle \alpha }}
A|_{x_{t_\alpha }=x_{t_{\alpha -}}\wedge x_{t_{\alpha +}}}\right)
+ \left( \sum_{\beta \in{\displaystyle
\beta }}A|_{x_{t^*_\beta
}=x_{t_{\beta -}}\wedge x_{t^*_{\beta +}}-x_{t^*_{\beta -}}\wedge x_{t_{\beta
+}}- [x_{t_{\beta -}},x_{t_{\beta +}}]}\right) -
(x_{t_-}-x_{t_+})\wedge A.
\eqno(\numb\label{eq825})
$$
$\square$
\end{remark}

\section{Differential bialgebra of bracket star-diagrams}\label{s9}

In this section we define the structure of differential bialgebras
on the complexes of bracket star-diagrams. We conjecture that this structure
is compatible with the corresponding bialgebra structure on the
homology space of the long knots space; the corresponding conjectures are 
given in~\cite{T1,T2,T3}.  

The cases of odd and even $d$
will be considered simultaneously.

Let $D$ be a bracket star-diagram, $T$ be a real number. Define a digram
$D^T$ as the diagram obtained from $D$ by the translation of $\R^1$ ($D^T$ 
is equal to $D$):
$$
t\mapsto t+T.
\eqno(\numb\label{eq91})
$$
Let now $A$ and $B$ be two bracket star-diagrams, we define their product
$A*B$ as the diagram $A^{-T}\cdot B^T$ in the case of odd $d$, and as 
the diagram $A^{-T}\wedge B^T$ in the case of even $d$, $T$ being a very large
positive number. This product resembles the product in the space of long 
knots:

\hspace*{25mm}
\PSbox{pic1.pstex}{10mm}{15mm} 
\begin{picture}(10,10)
\put(-11,12){$A$}
\put(52,12){$B$}
\put(134,12){$A$}
\put(197,12){$B$}
\end{picture}
\bigskip

\nopagebreak

\centerline{(Figure~\numb\label{f92})}

\bigskip

It follows from Theorem~\ref{t34} and Corollary~\ref{cor56} that the 
homology algebra of the differential algebra of bracket star-diagrams is 
commutative in the case of even $d$ for any commutative ring of coefficients.
This is also true over $\Q$ in the case of odd $d$, see Theorem~\ref{t510}.

\begin{conjecture}\label{con82}
Over $\Z$ in the case of odd $d$, the homology algebra of the differential
algebra of the bracket star-diagrams is not commutative.~$\square$
\end{conjecture}

Now we will define a comultiplication on the space of bracket star-diagrams.
In the case of odd $d$ the coproduct $\Delta$ of any
diagram $A=A_1\cdot A_2\ldots A_k$, where $A_i$, $1\le i\le k$, are brackets,
is defined as follows:
$$
\Delta A=\Delta(A_1\cdot A_2\ldots A_k):=
\sum\limits_{{I\sqcup J=\{1,\ldots,k\}
\atop I=\{i_1<\ldots<i_l\}}\atop J=\{j_1<\ldots <j_{k-l}\}}
(-1)^\epsilon A_{i_1}\cdot\ldots\cdot A_{i_l}\otimes
A_{j_1}\cdot\ldots\cdot A_{j_{k-l}},
\eqno(\numb\label{eq93})
$$
Where $\epsilon=\sum_{i_p>j_q}\tilde A_{i_p}\tilde A_{j_q}$.
  
In the case of even $d$
$$
\Delta (A)=\Delta(A_1\wedge\ldots\wedge A_k):=
\sum\limits_{{I\sqcup J=\{1,\ldots,k\}
\atop I=\{i_1<\ldots<i_l\}}\atop J=\{j_1<\ldots <j_{k-l}\}}
(-1)^\epsilon A_{i_1}\wedge\ldots\wedge A_{i_l}\otimes
A_{j_1}\wedge\ldots\wedge A_{j_{k-l}},
\eqno(\numb\label{eq94})
$$
where $\epsilon=\sum_{i_p>j_q}(\tilde A_{i_p}+1)(\tilde A_{j_q}+1)$.

In other words, $\Delta$ is the standard symmetric coalgebra coproduct.

It can be easily verified that the operations $\Delta$, $*$, $\partial$ define
a differential bialgebra structure on the space of bracket star-diagrams.

\vspace{2mm}

\hspace{3cm} Victor Tourtchine

\hspace{3cm} Independent University of Moscow,

\hspace{3cm} University of Paris 7

\hspace{3cm} Russia, 121002 Moscow,

\hspace{3cm} B.Vlassjevskij 11, MCCME

\hspace{3cm} e-mail: turchin@mccme.ru, tourtchi@acacia.ens.fr

\end{document}